\documentclass[12pt]{amsart}
\numberwithin{equation}{section}

\usepackage{amssymb}
\def\cb{{\mathcal B}}

\def\ce{{\mathcal E}}

\def\ch{{\mathcal H}}

\def\ck{{\mathcal K}}

\def\cs{{\mathcal S}}

\def\ga{{\mathfrak A}}

\def\gam{{\mathfrak M}}
\def\gn{{\mathfrak N}}

\def\bbf{{\mathbb F}}

\def\bt{{\mathbb T}}

\def\a{\alpha}

\def\g{\gamma} \def\G{\Gamma}
  \def\D{\Delta}

\def\eps{\varepsilon}
\def\l{\lambda} 

\def\m{\mu}

\def\r{\rho}
\def\s{\sigma} \def\S{\Sigma}
\def\t{\tau}
\def\f{\varphi} 
  
\def\om{\omega} \def\Om{\Omega}

\def\id{\hbox{id}}

\newtheorem{thm}{Theorem}[section]
\newtheorem{lem}[thm]{Lemma}

\newtheorem{prop}[thm]{Proposition}

\def\ad{\mathop{\rm Ad}}
\def\aut{\mathop{\rm Aut}}

\def\di{\mathop{\rm d}\!}

\def\oots{\overline{\otimes}}
\newcommand{\ots}[1]{\otimes_{\mathop{\rm #1}}}
\newcommand{\nn}{\nonumber}

\def\tr{\mathop{\rm Tr}}

\begin{document}

\title[ergodic theory]
{new topics in ergodic theory}
\author{Francesco Fidaleo}
\address{Francesco Fidaleo,
Dipartimento di Matematica,
Universit\`{a} di Roma Tor Vergata, 
Via della Ricerca Scientifica 1, Roma 00133, Italy} \email{{\tt
fidaleo@mat.uniroma2.it}}


\begin{abstract}
The entangled ergodic theorem concerns the study of the convergence in the strong, 
or merely weak operator topology, of the multiple Cesaro mean
$$
\frac{1}{N^{k}}\sum_{n_{1},\dots,n_{k}=0}^{N-1}
U^{n_{\a(1)}}A_{1}U^{n_{\a(2)}}\cdots 
U^{n_{\a(2k-1)}}A_{2k-1}U^{n_{\a(2k)}}\,,
$$    
where $U$ is a unitary operator acting on the Hilbert space $\ch$, 
$\a:\{1,\dots, m\}\mapsto\{1,\dots, k\}$ is a partition of the set 
made of $m$ elements in $k$ parts, and 
finally $A_{1},\dots,A_{2k-1}$ are bounded operators acting on $\ch$.

While reviewing recent results about the entangled ergodic theorem, 
we provide some natural applications to dynamical systems based 
on compact operators. 

Namely,
let $(\ga,\a)$ be a $C^{*}$--dynamical 
system, where $\ga=\ck(\ch)$, and $\a=\ad{}\!_{U}$ is an automorphism implemented 
by the unitary $U$. 

We show that
$$
\lim_{N\to+\infty}\frac{1}{N}\sum_{n=0}^{N-1}\a^{n}=E\,,
$$
pointwise in the weak topology of $\ck(\ch)$. Here, $E$ is a conditional 
expectation projecting onto the $C^{*}$--subalgebra 
$$
\bigg(\bigoplus_{z\in\s_{\mathop{\rm pp}}(U)}
E_{z}\cb(\ch)E_{z}\bigg)\bigcap\ck(\ch)\,.
$$

If in addition $U$ is weakly mixing with $\Om\in\ch$ 
the unique up to a phase, invariant vector under $U$ and 
$\om=\langle\cdot\,\Om,\Om\rangle$, we have the following recurrence 
result. If $A\in\ck(\ch)$ fulfils
$\om(A)>0$, and $0<m_{1}<m_{2}<\cdots<m_{l}$ are 
natural numbers kept fixed, then there exists an $N_{0}$ such that 
$$
\frac{1}{N}\sum_{n=0}^{N-1}\om(A\a^{nm_{1}}(A)\a^{nm_{2}}(A)\cdots
\a^{nm_{l}}(A))>0
$$
for each $N>N_{0}$.
\end{abstract}

\maketitle

\section{introduction}

Recently, it was shown that same ergodic properties of classical 
dynamical systems fail to be true by passing to noncommutative setting.
It is then of interest to understand among the various 
ergodic properties, which ones survive by passing from the classical 
to the quantum case. We mention the pivotal paper \cite{NSZ}, where such 
an investigation is 
carried out for some basic recurrence, as well as multiple mixing properties.

Notice that 
it is in general unclear what should be the right quantum counterpart 
of a classical ergodic property. As an example, we mention   
the property of the 
convergence to the equilibrium (i.e. ergodicity for an invariant 
state $\om$)
$$
\lim_{N\to+\infty}\frac{1}{N}\sum_{n=0}^{N-1}\om\left(B^{*}\a^{n}(A)B\right)
=\om\left(B^{*}B\right)\om(A)
$$
suggested by the quantum physics, and the standard definition of 
ergodicity
$$
\lim_{N\to+\infty}\frac{1}{N}\sum_{n=0}^{N-1}\om(A\a^{n}(B))
=\om(A)\om(B)\,.
$$
We refer the reader to \cite{FL}, Proposition 1.1 for further details.

A notion which is meaningful in quantum setting is that of 
entangled ergodic theorem, introduced in \cite{AHO} in 
connection with the central limit theorem for suitable sequences of elements of the group 
$C^{*}$--algebra of the free group $\bbf_{\infty}$ on infinitely many generators. 

The entangled ergodic theorem can be clearly formulated in the 
following way. Let $U$ be 
a unitary operator acting on the Hilbert space $\ch$, and for $m\geq k$, 
$\a:\{1,\dots,m\}\mapsto\{1,\dots, k\}$ a partition of the 
set $\{1,\dots,m\}$ in $k$ parts. 
The entangled ergodic theorem concerns 
the convergence in the strong, or merely weak operator topology, 
of the multiple Cesaro mean
\begin{equation}
\label{0}
\frac{1}{N^{k}}\sum_{n_{1},\dots,n_{k}=0}^{N-1}
U^{n_{\a(1)}}A_{1}U^{n_{\a(2)}}\cdots 
U^{n_{\a(m-1)}}A_{m-1}U^{n_{\a(m)}}\,,
\end{equation}
$A_{1},\dots,A_{m-1}$ being bounded operators acting on $\ch$.

Notice that expressions like \eqref{0} naturally appear also in \cite{NSZ}
relatively to the study of the behaviour of the 
multiple correlations. Just by considering the simplest case 
of the  
partition of the empty set, the limit of the 
Cesaro mean in \eqref{0} reduces itself to the well--known mean ergodic theorem 
due to John von Neumann (cf. \cite{RS}) 
\begin{equation}
\label{jvn}
\mathop{\rm s\!-\!lim}_{N\to+\infty}\frac{1}{N}\sum_{n=0}^{N-1}U^{n}=E_{1}\,,
\end{equation}
$E_{1}$ being the selfadjoint projection onto the eigenspace of the 
invariant vectors for $U$.

Some applications of the entangled ergodic theorem are discussed below. 
Apart from the other potential applications to the study of the ergodic 
properties of quantum dynamical systems, the entangled ergodic theorem 
is a fascinating self--contained mathematical problem. It is certainly true if the 
spectrum $\s(U)$ of $U$ is finite. Some very special cases for which 
it holds true are listed in \cite{L}. 

The first part of the present paper, based on \cite{F, F1}, is devoted to review the known 
results on the entangled ergodic theorem. 

We start by considering 
the sufficiently general situation when the operators 
$A_{1},\dots,A_{m-1}$ in \eqref{0} are compact (cf. \cite{F}).

Then we pass to the case when the unitary $U$ is 
almost periodic (i.e. $\ch$ is generated by the eigenvectors of 
$U$), and $\a:\{1,\dots,2k\}\mapsto\{1,\dots, k\}$ a pair--partition, 
without any condition on the operators $A_{1},\dots,A_{2k-1}$ (cf. 
\cite{F1}).

Another interesting situation arises from the generalization to the 
noncommutative setting, of the ergodic theorem of H. Furstenberg relative 
to diagonal measure (cf. \cite{F1, Fu1}). By using such a result, 
we can treat the followig situation. Let 
$M$ be a von Neumann algebra equipped with the adjoint action of an ergodic unitary $U$,
and a standard vector $\Om$ which is invariant under $U$. Let $M'$ be the 
commutant von Neumann algebra of $M$. 
In this situation, the Cesaro mean
\begin{equation}
\label{00}
\frac{1}{N}\sum_{n=0}^{N-1}U^{n}AU^{n}
\end{equation}
converges in the strong operator topology for each $A\in M\bigcup M'$ 
(cf. \cite{F1}). Notice that \eqref{00} is the particular case of \eqref{0} relative to 
the trivial 
pair--partition of two elements.

The second part of the present paper concerns the application of the 
entangled ergodic theorem, as well as some lines of its proof, to the 
investigation of ergodic properties of $C^{*}$--dynamical systems 
based on compact operators. More precisely, let $(\ga,\a)$ be a $C^{*}$--dynamical 
system, where $\ga=\ck(\ch)$, and $\a=\ad{}\!_{U}$ is an automorphism implemented 
by the unitary $U$. We show that
$$
\lim_{N\to+\infty}\frac{1}{N}\sum_{n=0}^{N-1}\a^{n}=E\,,
$$
pointwise in the weak topology of $\ck(\ch)$. Here, $E$ is a conditional 
expectation projecting onto 
the $C^{*}$--subalgebra 
$$
\bigg(\bigoplus_{z\in\s_{\mathop{\rm pp}}(U)}
E_{z}\cb(\ch)E_{z}\bigg)\bigcap\ck(\ch)\,.
$$

If in addition $U$ is weakly mixing with 
$\Om\in\ch$ 
the unique up to a phase, invariant vector under $U$, we can 
consider the weakly mixing $C^{*}$--dynamical 
system $(\ga,\a,\om)$ where
$\om=\langle\cdot\,\Om,\Om\rangle$. We prove the following recurrence 
result. If $A\in\ck(\ch)$ satisfies
$\om(A)>0$, and $0<m_{1}<m_{2}<\cdots<m_{l}$ are 
natural numbers kept fixed, then there exists an $N_{0}$ such that 
$$
\frac{1}{N}\sum_{n=0}^{N-1}\om(A\a^{nm_{1}}(A)\a^{nm_{2}}(A)\cdots
\a^{nm_{l}}(A))>0
$$
for each $N>N_{0}$.

We end the present section with some notations and definitions useful in 
the sequel.

The convergence in the weak, respectively strong operator topology 
(see e.g. \cite{SZ, T}) of 
a net $\{A_{\a}\}_{\a\in J}\subset\cb(\ch)$ is denoted respectively as
$$
\mathop{\rm w\!-\!lim}_{\a}A_{\a}=A\,,\quad
\mathop{\rm s\!-\!lim}_{\a}A_{\a}=A\,.
$$

Let $U$ be a unitary operator acting on $\ch$. Consider the resolution of the identity  
$\{E(B)\,:\,B\,\,\text{Borel subset of}\,\,\bt\}$ of $U$ (cf. 
\cite{TL}, Section VII.7). Denote with an abuse of notation, 
$E_{z}:=E(\{z\})$. Namely, $E_{z}$ is nothing but the selfadjoint 
projection on the eigenspace corresponding to the eigenvalue 
$z$ in the unit circle $\bt$. Denote $\s_{\mathop{\rm pp}}(U):=
\big\{z\in\bt\,:\,z\,\text{is an eigenvalue of}\,U\big\}$ (cf. 
\cite{RS}).

The unitary $U$ is said to be {\it ergodic} if the 
fixed--point subspace $E_{1}\ch$ is one 
dimensional. By the mean ergodic theorem \eqref{jvn}, it is equivalent to the 
existence of a unit vector $\xi_{0}\in\ch$ such that 
$$
\lim_{N\to+\infty}\frac{1}{N}\sum_{n=0}^{N-1}U^{n}\xi=
\langle\xi,\xi_{0}\rangle\xi_{0}\,,
$$
or equivalently,
$$
\lim_{N\to+\infty}\frac{1}{N}\sum_{n=0}^{N-1}\langle U^{n}\xi,\eta\rangle=
\langle\xi,\xi_{0}\rangle\langle\xi_{0},\eta\rangle\,.
$$

The unitary $U$ is said to be {\it weakly mixing} if there exists 
a unit vector $\xi_{0}\in\ch$ such that 
$$
\lim_{N\to+\infty}\frac{1}{N}\sum_{n=0}^{N-1}\big|\langle U^{n}\xi,\eta\rangle
-\langle\xi,\xi_{0}\rangle\langle\xi_{0},\eta\rangle\big|=0\,.
$$

A unitary $U$ is weakly mixing if and only if 
$\s_{\mathop{\rm pp}}(U)=\{1\}$ and 
$E_{1}=\langle\,\cdot\,,\xi_{0}\rangle\xi_{0}$, see e.g. \cite{NSZ}, 
Proposition 5.4.

The unitary $U$ is said to be {\it almost periodic} if
$\ch=\ch_{\mathop{\rm ap}}^{U}$, $\ch_{\mathop{\rm ap}}^{U}$ being 
the closed subspace 
consisting of the vectors having relatively norm--compact orbit 
under $U$. It is seen in \cite{NSZ} that $U$ is almost periodic if and 
only if $\ch$ is generated by the eigenvectors of $U$.

For a (discrete) $C^*$--dynamical system we mean a pair $\big(\ga,\a,\big)$ 
consisting of a $C^*$-algebra $\ga$, and an automorphism $\a$ of $\ga$. 
If in addition, a state $\om\in\cs(\ga)$ invariant under the action of 
$\a$ is kept fixed, we consider also $C^*$--dynamical systems 
consisting of a triplet $\big(\ga,\a,\om\big)$. 

A $C^*$--dynamical system $\big(\ga,\a,\om\big)$
is said to be {\it ergodic} if for each 
$A,B\in\ga$,
$$
\lim_{N\to+\infty}\frac{1}{N}\sum_{n=0}^{N-1}\om(A\a^{n}(B))=
\om(A)\om(B)\,.
$$
It is said to be {\it weakly mixing} if
$$
\lim_{N\to+\infty}\frac{1}{N}\sum_{n=0}^{N-1}\big|\om(A\a^{n}(B))-
\om(A)\om(B)\big|=0
$$
for each $A,B\in\ga$. 

Let $\big(\ch,\pi,U,\Om\big)$ be the GNS covariant representation 
(cf. \cite{T}, Section I.9)  
canonically associated to the dynamical system under 
consideration. Then $\big(\ga,\a,\om\big)$ is ergodic (respectively 
weakly mixing) if and only if $U$ is ergodic (respectively 
weakly mixing), see e.g. \cite{NSZ}.

\section{the entangled ergodic theorem}

The present section, based on \cite{F, F1}, is devoted to review the known 
results on the entangled ergodic theorem. 

\subsection{case of compact operators}${}$\\

We start with the entangled ergodic theorem 
for general partitions of any finite set $\{1,\dots,m\}$, and for 
compact operators $\{A_{1},\dots,A_{m-1}\}$.

Let $U\in\cb(\ch)$ be a unitary operator, and for $m\geq k$, 
$\a:\{1,\dots, m\}\mapsto\{1,\dots, k\}$ a partition of the 
set $\{1,\dots, m\}$ in $k$ parts.\footnote{A partition 
$\a:\{1,\dots,m\}\mapsto\{1,\dots, k\}$
of the set made of $m$ elements in $k$ parts is nothing but a surjective 
map, the parts of $\{1,\dots,m\}$ being the preimages 
$\{\a^{-1}(\{j\})\}_{j=1}^{k}$.}
\begin{thm}(cf. \cite{F}, Theorem 2.6)\\ 
\label{main1}
For $m\geq k$, let $\a:\{1,\dots,m\}\mapsto\{1,\dots,k\}$ be a partition of the 
set $\{1,\dots,m\}$. If $\{A_{1},\dots,A_{m-1}\}\subset\ck(\ch)$, then the
ergodic average 
\begin{equation*}
\frac{1}{N^{k}}\sum_{n_{1},\dots,n_{k}=0}^{N-1}
U^{n_{\a(1)}}A_{1}U^{n_{\a(2)}}\cdots 
U^{n_{\a(m-1)}}A_{m-1}U^{n_{\a(m)}}
\end{equation*}
converges in the weak operator topology to some bounded operator 
$S_{\a;A_{1},\dots,A_{m-1}}\in\cb(\ch)$.
\end{thm}
\begin{proof}
Define
$$
\G_{N}:=\frac{1}{N^{k}}\sum_{n_{1},\dots,n_{k}=0}^{N-1}
U^{n_{\a(1)}}A_{1}U^{n_{\a(2)}}\cdots 
U^{n_{\a(m-1)}}A_{m-1}U^{n_{\a(m)}}\,.
$$
Notice that
$$
\|\G_{N}\|\leq\prod_{j=1}^{m-1}\|A_{j}\|\,.
$$

By Theorem II.1.3 of \cite{T}, it is then enough to show that the 
$\langle\G_{N}x,y\rangle$ converges 
for each fixed $x,y\in\ch$. On the other hand, we can approximate 
the $A_{j}$ by finite rank operators. Namely, put $K:=\max_{1\leq j\leq m-1}\|A_{j}\|$. 
Choose finite rank 
operators $A^{\eps}_{j}$, such that $\|A^{\eps}_{j}\|\leq K$ and
$\|A_{j}-A^{\eps}_{j}\|<{\displaystyle\frac{\eps}{4(m-1)K^{m-2}\|x\|\|y\|}}$,
$j=1,\dots,m-1$. We have with obvious notations
\begin{align*}
&|\langle\G_{N}x,y\rangle-\langle\G_{M}x,y\rangle|\leq
\|\G_{N}-\G^{\eps}_{N}\|+\|\G_{M}-\G^{\eps}_{M}\|\\
+&|\langle\G^{\eps}_{N}x,y\rangle-\langle\G^{\eps}_{M}x,y\rangle|
\leq\frac{\eps}{2}+|\langle\G^{\eps}_{N}x,y\rangle-\langle\G^{\eps}_{M}x,y\rangle|\,.
\end{align*}
Thus, it is enough to show that  
$$
\bigg\langle\frac{1}{N^{k}}\sum_{n_{1},\dots,n_{k}=0}^{N-1}
U^{n_{\a(1)}}A_{1}U^{n_{\a(2)}}\cdots 
U^{n_{\a(m-1)}}A_{m-1}U^{n_{\a(m)}}x,y\bigg\rangle
$$
converges for every $x,y\in\ch$, whenever the $A_{j}$ are rank one
operators. By using the explicit computations in \cite{L}, we obtain in this situation,
\begin{align*}
&\bigg\langle\frac{1}{N^{k}}\sum_{n_{1},\dots,n_{k}=0}^{N-1}
U^{n_{\a(1)}}A_{1}U^{n_{\a(2)}}\cdots 
U^{n_{\a(m-1)}}A_{m-1}U^{n_{\a(m)}}x,y\bigg\rangle\\
=&\prod_{j=1}^{k}\frac{1}{N}\sum_{n_{j}=0}^{N-1}
\prod_{\{p\,|\,\a(p)=j\}}\big\langle U^{n_{j}}x_{p,j},y_{p,j}\big\rangle\\
=&\prod_{j=1}^{k}\iint\cdot\cdot\int_{\bt^{|\a^{-1}\{j\}|}}
\bigg(\frac{1}{N}\sum_{n_{j}=0}^{N-1}
\big(\prod_{\{p\,|\,\a(p)=j\}}z_{p}\big)^{n_{j}}\bigg)
\prod_{\{p\,|\,\a(p)=j\}}\langle E(\di z_{p})x_{p,j},y_{p,j}\rangle\\
&\longrightarrow_{{}_{N}}
\prod_{j=1}^{k}\iint\cdot\cdot\int_{\bt^{|\a^{-1}\{j\}|}}
\chi_{\{1\}}\bigg(\prod_{\{p\,|\,\a(p)=j\}}z_{p}\bigg)
\prod_{\{p\,|\,\a(p)=j\}}\langle E(\di z_{p})x_{p,j},y_{p,j}\rangle
\end{align*}
where we have 
used the Lebesgue dominated convergence theorem. 
Here, the $x_{p,j}$, $y_{p,j}$ are vectors uniquely determined by the 
rank one operators $A_{1},\dots,A_{m-1}$ and vectors $x,y$, 
and finally $\chi_{\D}$ denotes the indicator of the set $\D$.
\end{proof}

It was shown in \cite{F} that if $\a:\{1,\dots,2k\}\mapsto\{1,\dots,k\}$ 
is a pair--partition, we can explicitely write the formula for 
$S_{\a;A_{1},\dots,A_{2k-1}}\in\cb(\ch)$. Namely, define
\begin{equation*}
\s_{\mathop{\rm pp}}^{\mathop{\rm a}}(U):=
\big\{z\in\s_{\mathop{\rm pp}}(U)\,:\,zw=1\,\text{for 
some}\,w\in\s_{\mathop{\rm pp}}(U)\big\}
\end{equation*}
Then we have
\begin{align}
\label{trc}
\mathop{\rm w\!-\!lim}_{N}&\bigg\{\frac{1}{N^{k}}\sum_{n_{1},\dots,n_{k}=0}^{N-1}
U^{n_{\a(1)}}A_{1}U^{n_{\a(2)}}\cdots 
U^{n_{\a(2k-1)}}A_{2k-1}U^{n_{\a(2k)}}\bigg\}\nn\\
=&\sum_{z_{1},\dots,z_{k}\in\s_{\mathop{\rm pp}}^{\mathop{\rm a}}(U)}
E_{z^{\#}_{\a(1)}}A_{1}E_{z^{\#}_{\a(2)}}\cdots
E_{z^{\#}_{\a(2k-1)}}A_{2k-1}E_{z^{\#}_{\a(2k)}}\,,
\end{align} 

Here, the pairs $z^{\#}_{\a(i)}$ are alternatively $z_{j}$ and
$\bar z_{j}$ whenever $\a(i)=j$,
$E_{z}$ is the selfadjoint 
projection on the eigenspace corresponding to the eigenvalue 
$z\in\s_{\mathop{\rm pp}}(U)$,\footnote{If for example, $\a$ is the pair--partition
$\{1,2,1,2\}$ of four elements, 
$$
S_{\a;A,B,C}=\sum_{z,w\in\s_{\mathop{\rm pp}}^{\mathop{\rm a}}(U)}
E_{z}AE_{w}BE_{\bar z}CE_{\bar w}\,.
$$} 
and finally
the sum in the r.h.s. is understood as the limit in the weak 
operator topology of the net obtained by considering the finite 
truncations of the r.h.s. of \eqref{trc} (cf. \cite{F}, Proposition 
2.3). Notice that \eqref{trc} cannot be extended to the whole 
$\cb(\ch)$, see the example in pag. 8 of \cite{NSZ}.

\subsection{almost periodic case}${}$\\

Another case for which the entangled ergodic theorem can be proved is 
the almost periodic case, that is when the Hilbert space is generated 
by the eigenvectors of the unitary $U$. In this situation, we have no 
conditions on the bounded operators appearing in \eqref{0}.
\begin{thm}(cf. \cite{F1}, Theorem 2.6)\\
\label{qper2}    
Suppose that the 
dynamics induced by the unitary $U$ on $\ch$ is almost periodic. Then 
for each $A_{1},\dots,A_{2k-1}\in\cb(\ch)$,
\begin{align*}
\mathop{\rm s\!-\!lim}_{N}&\bigg\{\frac{1}{N^{k}}\sum_{n_{1},\dots,n_{k}=0}^{N-1}
U^{n_{\a(1)}}A_{1}U^{n_{\a(2)}}\cdots 
U^{n_{\a(2k-1)}}A_{2k-1}U^{n_{\a(2k)}}\bigg\}\\
=&S_{\a;A_{1},\dots,A_{2k-1}}\,.
\end{align*} 
\end{thm}
\begin{proof}
To simplify, we treat the case of the partition $\{1,2,1,3,2,3\}$, the general 
case follows the same lines of this case.
Fix $\eps>0$, and suppose that $A,B,C,D,F\in\cb(\ch)$ have norm  
one. Let $I_{\eps}$ be such that 
$$
\bigg\|x-\sum_{\s\in I_{\eps}}E_{\s}x\bigg\|<\eps\,.
$$
For each $\s\in I_{\eps}$, let $I_{\eps}(\s)$ be such that 
$$
\bigg\|FE_{\s}x-\sum_{\t\in 
I_{\eps}(\s)}E_{\t}FE_{\s}x\bigg\|<\frac{\eps}{|I_{\eps}|}\,.
$$
Choose $N_{\eps}$ such that 
$$
\bigg\|\bigg(\frac{1}{N}\sum_{n=0}^{N-1}(\s U)^{n}-E_{\bar\s}\bigg)
DE_{\t}FE_{\s}x\bigg\|<\frac{\eps}{{\displaystyle\sum_{\s\in I_{\eps}}|I_{\eps}(\s)|}}\,,
$$
whenever $N>N_{\eps}$ and $\s\in I_{\eps}$, $\t\in I_{\eps}(\s)$.  
For each $\s\in I_{\eps}$, $\t\in I_{\eps}(\s)$, let $I_{\eps}(\s,\t)$
be such that 
$$
\bigg\|CE_{\bar\s}DE_{\t}FE_{\s}x-\sum_{\r\in I_{\eps}(\s,\t)}E_{\r}C
E_{\bar\s}DE_{\t}FE_{\s}x\bigg\|<\frac{\eps}
{{\displaystyle\sum_{\s\in I_{\eps}}|I_{\eps}(\s)|}}\,.
$$
Then
\begin{align*}
&\bigg\|\frac{1}{N^{3}}\sum_{k,m,n=0}^{N-1}U^{k}AU^{m}BU^{k}CU^{n}DU^{m}FU^{n}x
-S_{\a;A,B,C,D,F}x\bigg\|\\
\leq&5\eps
+\sum_{\s\in I_{\eps}}\sum_{\t\in I_{\eps}(\s)}\sum_{\r\in I_{\eps}(\s,\t)}
\bigg\|\bigg(\frac{1}{N}\sum_{k=0}^{N-1}(\r U)^{k}\bigg)A
\bigg(\frac{1}{N}\sum_{m=0}^{N-1}(\t U)^{m}\bigg)\\
\times&BE_{\r}C
E_{\bar\s}DE_{\t}FE_{\s}x-E_{\bar\r}AE_{\bar\t}B
E_{\r}CE_{\bar\s}DE_{\t}FE_{\s}x\bigg\|\,.
\end{align*}
Taking the limsup on both sides, we obtain the assertion by the mean ergodic theorem
\eqref{jvn}.
\end{proof}
 
\subsection{diagonal measures}${}$\\

We treat the natural generalization to the 
quantum case of the celebrated result due to H. Furstenberg relative to 
the diagonal measures (cf. \cite{Fu1}, Section 
4.4).

We start with a $C^*$--dynamical system $\big(\ga,\a,\om\big)$, 
together with its GNS covariant representation $\big(\ch,\pi,U,\Om\big)$. 
Denote $M:=\pi(\ga)''$,  
the von Neumann algebra acting on $\ch$ generated by 
the representation $\pi$. The commutant von Neumann algebra is denoted as $M'$.
Suppose further that the support $s(\om)$ in $\ga^{**}$
is central. The last property simply means that $\Om$ is 
separating for $\pi(\ga)''$, see e.g. \cite{SZ}, Section 10.17. 

Let $\gam:=M\ots{max}M'$ be the completion of the algebraic 
tensor product $\gn:=M\otimes M'$ w.r.t. the maximal $C^*$--norm (cf. 
\cite{T}, Section IV.4). It is easily seen that on $\gam$ the 
following two states are automatically well--defined. The first one is 
the canonical product state 
$$
\f(A\otimes B):=\langle A\Om,\Om\rangle\langle B\Om,\Om\rangle\,,
\quad A\in M\,,\,\,B\in M'\,.
$$
The second one is uniquely defined by
$$
\psi(A\otimes B):=\langle AB\Om,\Om\rangle\,,
\quad A\in M\,,\,\,B\in M'\,.
$$

The state $\psi$ can be considered the (quantum analogue of the) ``diagonal 
measure'' of the ``measure'' $\f$.

On $\gam$ is also uniquely defined the automorphism
$$
\g:=\ad{}\!_{U}\otimes\ad{}\!_{U^{2}}\,,
$$
see \cite{T}, Proposition IV.4.7. Of course, $\big(\gam,\g,\f\big)$ is 
a $C^*$--dynamical system whose GNS covariant representation is 
precisely 
$\big(\ch\otimes\ch,\id\otimes\id,U\otimes U^{2},\Om\otimes\Om\big)$.
Denote $E_{1}$ the selfadjoint 
projection onto the invariant vectors under $U\otimes U^{2}$.
Notice that the 
$*$--subalgebra $\gn$ is globally stable under the action of $\g$.

In addition, again by Proposition IV.4.7 of \cite{T}, 
$$
\s(A\otimes B):=AB\,,
\quad A\in M\,,\,\,B\in M'\,.
$$
uniquely defines a representation of $\gam$ on $\ch$ such that 
$\big(\ch,\s,\Om\big)$ is precisely the GNS representation of the 
state $\psi$.

Let $A\in M$, $B\in M'$. Then by the mean ergodic theorem \eqref{jvn},
\begin{align*}
\frac{1}{N}\sum_{n=0}^{N-1}\psi(\g^{n}(A&\otimes B))
=\frac{1}{N}\sum_{n=0}^{N-1}\langle AU^{n}B\Om,\Om\rangle\\
&\equiv\bigg\langle A\bigg(\frac{1}{N}\sum_{n=0}^{N-1}U^{n}\bigg)B\Om,
\Om\bigg\rangle\\
&\longrightarrow\langle A\Om,\Om\rangle\langle B\Om,\Om\rangle
\equiv\f(A\otimes B)\,.
\end{align*}

According with Definition 4.1 of \cite{F1} (see also \cite{Fu1}, Definition 4.4), 
this means that the state 
$\psi\in\cs(\gam)$ is 
generic for $\big(\gam,\g,\f\big)$ w.r.t. $\gn$. In addition, define 
$$
\S:=\{(z,w)\in\s_{\mathop{\rm pp}}
(U)\times\s_{\mathop{\rm pp}}(U)\,:\,
zw^{2}=1\}\,.
$$
Then by Lemma 4.18 of \cite{Fu1},
$$
E_{1}=\bigoplus_{s\in\S}E^{U}_{z_{s}}\otimes E^{U}_{w_{s}}\,,
$$
$E^{U}_{z}$ being the selfadjoint projection onto 
the eigenspace of $U$ corresponding to the eigenvalue $z$. As $U$ is 
supposed to be 
ergodic, by Proposition 2.2 of \cite{F1}, $E^{U}_{z}\ch$ is one dimensional, 
and $E^{U}_{z_{s}}\ch$ and $E^{U}_{w_{s}}\ch$ are generated by 
$V_{z_{s}}\Om$, $W_{w_{s}}\Om$, where $V_{z_{s}}$ and $W_{w_{s}}$ are 
unitaries of $M_{z_{s}}:=\{A\in M\,:\,UAU^{-1}=z_{s}A\}$, 
$(M')_{w_{s}}:=\{B\in M'\,:\,UBU^{-1}=w_{s}B\}$ respectively. Thus,
$E^{U}_{z_{s}}\ch\oots E^{U}_{w_{s}}\ch$ is one dimensional, and it 
is generated by $V_{z_{s}}\Om\otimes W_{w_{s}}\Om$. This means that 
$\gn\Om\bigcap E_{1}\ch\oots\ch$ is dense in $E_{1}\ch\oots\ch$. 
Then the map 
$$
\sum_{j}A_{j}\Om\otimes B_{j}\Om\in\gn\Om\bigcap E_{1}\ch\oots\ch
\mapsto\sum_{j}A_{j}B_{j}\Om\in\ch
$$
uniquely defines a partial isometry $V:\ch\oots\ch\mapsto\ch$ such 
that $V^{*}V=E_{1}\ch\oots\ch$. For $\xi,\eta\in\ch$, such an isometry has the form
\begin{equation}
\label{isom}
V\left(\xi\otimes\eta\right)=\sum_{\{z,w\in\s(U)\,:\,zw^{2}=1\}}
\langle \xi,V_{z}\Om\rangle\langle\eta,W_{w}\Om\rangle V_{z}W_{w}\Om\,,
\end{equation}
where $V_{z}\Om$, $V_{z}$ unitary of $M_{z}$
(equivalently $W_{z}\Om$, $W_{z}$ unitary of $(M')_{z}$) generates the one 
dimensional subspace $E^{U}_{z}\ch$ for $z\in\s_{\mathop{\rm pp}}(U)$.
\begin{thm}(cf. \cite{F1}, Theorem 5.2)\\
\label{diam1}    
Let $\big(\ga,\a,\om\big)$ be an ergodic $C^*$--dynamical system such 
that its support $c(\om)$ in $\ga^{**}$ is central. Then for each 
$A\in M\bigcup M'$,
\begin{equation*}
\mathop{\rm s\!-\!lim}_{N\to+\infty}
\frac{1}{N}\sum_{n=0}^{N-1}U^{n}AU^{n}
=V\left(A\Om\otimes\,\cdot\,\right)\,,
\end{equation*}
where $V$ is the isometry given in \eqref{isom}.
\end{thm}
\begin{proof}
As $\psi$ is 
generic for $\big(\gam,\g,\f\big)$ w.r.t. $M\otimes M'$ and the 
last $*$--algebra is left stable by 
$\ad{}\!_{U}\otimes\ad{}\!_{U^{2}}$, we can apply Theorem 4.5 of 
\cite{F1} (see also \cite{Fu1}, Theorem 4.14 for the Abelian case) 
obtaining for $X\in M$, $Y\in M'$
$$
\lim_{N\to+\infty}\frac{1}{N}\sum_{n=0}^{N-1}U^{n}XU^{n}Y\Om
=V\left(X\Om\otimes Y\Om\right)\,.
$$

If $A\in M$ the proof follows as $\Om$ is cyclic for $M'$. By 
exchanging the role between $M$ and $M'$, we obtain the result 
whenever $A\in M'$ 
\end{proof}

If $\big(\ga,\a,\om\big)$ is weakly mixing and 
$0<m_{1}<m_{2}$ natural numbers, we prove following the same lines of 
Theorem \ref{diam1}, but in a different way from 
\cite{NSZ}, that 
\begin{equation}
\label{sssta}
\mathop{\rm s\!-\!lim}_{N\to+\infty}
\frac{1}{N}\sum_{n=0}^{N-1}U^{nm_{1}}AU^{nm_{2}}
=\langle A\Om,\Om\rangle\langle\,\cdot\,,\Om\rangle\Om
\end{equation}
for each $A\in M\bigcup M'$.

\section{applications}

We start with the following recurrence result which is a direct 
consequence of Theorem 1.3 of \cite{NSZ}. By \eqref{sssta}, we then have an 
alternative proof of it.
\begin{prop}
\label{recc1}
Let $\big(\ga,\a,\om\big)$ be a weakly mixing $C^*$--dynamical system such 
that its support $c(\om)$ in $\ga^{**}$ is central, and 
$0<m_{1}<m_{2}$ natural numbers. Consider $A\in\ga$ such that 
$\om(A)>0$

Then there exists an $N_{0}$ such that 
$$
\frac{1}{N}\sum_{n=0}^{N-1}\om(A\a^{nm_{1}}(A)\a^{nm_{2}}(A))>0
$$
for each $N>N_{0}$.
\end{prop}
\begin{proof}
We have by \eqref{sssta}, 
\begin{align*}
&\lim_{N\to+\infty}\frac{1}{N}\sum_{n=0}^{N-1}
\om(A\a^{nm_{1}}(A)\a^{nm_{2}}(A))\\
\equiv&\lim_{N\to+\infty}
\frac{1}{N}\sum_{n=0}^{N-1}\big\langle AU^{nm_{1}}AU^{n(m_{2}-m_{1})}A
\Om,\Om\big\rangle\\
=&\om(A)^{3}>0\,.
\end{align*}
\end{proof}

Now we pass to some interesting applications concerning compact 
operators.
\begin{prop}
\label{recc2}
Let $U$ be a weakly mixing unitary acting on the Hilbert space $\ch$, 
$A_{1},\dots,A_{k-1}\in\ck(\ch)$, and finally $m_{1},\dots,m_{k}$ 
fixed nonnull natural numbers. Then
\begin{align}
\label{wwl}
\mathop{\rm w\!-\!lim}_{N\to+\infty}
\frac{1}{N}\sum_{n=0}^{N-1}
U^{nm_{1}}A_{1}&U^{nm_{2}}A_{2}\cdots 
U^{nm_{k-1}}A_{k-1}U^{nm_{k}}\nn\\
=&E_{1}A_{1}E_{1}A_{2}\cdots E_{1}A_{k-1}E_{1}\,.
\end{align}
\end{prop}
\begin{proof}
As $U$ is weakly mixing, $\s_{\mathop{\rm pp}}(U)=\{1\}$ and 
$E_{1}=\langle\,\cdot\,,\Om\rangle\Om$ for a unique up to a phase unit 
vector. In addition, we can approximate 
the $A_{j}$ by finite rank operators as explained in Theorem 
\ref{main1}. We now decompose $U$ as 
\begin{equation}
\label{wwwl}
U=\langle\,\cdot\,,\Om\rangle\Om+E_{1}^{\perp}U\,,
\end{equation}
where $E_{1}^{\perp}$ is the selfadjoint projection onto the closed subspace on 
which $U$ has purely continuous spectrum. By inserting \eqref{wwwl} in
\eqref{wwl}, we obtain an addendum containing in all place the piece
$\langle\,\cdot\,,\Om\rangle\Om$, the last coinciding with 
$E_{1}A_{1}E_{1}A_{2}\cdots E_{1}A_{k-1}E_{1}$. As we reduced the 
matter to the case when the $A_{j}$ are rank one operators, the remaining addenda 
contain a multiplicative factor of the form
\begin{equation}
\label{wwwll}
G_{N}:=\iint\cdot\cdot\int_{\bt^{j}}\left(\frac{1}{N}\sum_{n=0}^{N-1}
\big(z_{1}^{m_{1}}\cdots z_{j}^{m_{j}}\big)^{n}\right)\di\m_{1}(z_{1})\cdots
\di\m_{j}(z_{j})\,.
\end{equation}

In \eqref{wwwll} $1\leq j\leq k$ is fixed and depends on the addendum under 
consideration, and 
$$
\di\m_{l}(z_{l}):=\langle E(\di z_{l})x_{l},y_{l}\rangle\,,\quad 
1\leq l\leq j
$$
are bounded signed Borel measure without atoms. As
$$
\frac{1}{N}\sum_{n=0}^{N-1}(z_{1}^{m_{1}}\cdots z_{j}^{m_{j}})^{n}\longrightarrow
\chi_{\{1\}}(z_{1}^{m_{1}}\cdots z_{j}^{m_{j}})
$$
pointwise, by taking the limit 
in \eqref{wwwll}, we obtain by Lebesgue dominated convergence theorem 
and Fubini theorem,
$$
\lim_{N\to+\infty}G_{N}=\iint\cdot\cdot\int_{\bt^{j-1}}
f(z_{1},\dots,z_{j-1})\di\m_{1}(z_{1})\cdots\di\m_{j-1}(z_{j-1})
$$    
where 
$$
f(z_{1},\dots,z_{j-1}):=\m_{j}\big(\big\{z_{j}\,:\,
z_{1}^{m_{1}}\cdots z_{j}^{m_{j}}=1\big\}\big)\,.
$$

The proof follows as, for fixed $z_{1},\dots,z_{j-1}\in\bt$,  
$$
\#\big\{z_{j}\,:\,z_{1}^{m_{1}}\cdots z_{j}^{m_{j}}=1\big\}=m_{j}\,.
$$
\end{proof}

Notice that, if
$\s_{\mathop{\rm pp}}(U)=\emptyset$, then
$$
\mathop{\rm w\!-\!lim}_{N\to+\infty}
\frac{1}{N}\sum_{n=0}^{N-1}
U^{nm_{1}}A_{1}U^{nm_{2}}A_{2}\cdots 
U^{nm_{k-1}}A_{k-1}U^{nm_{k}}
=0\,.
$$   

The proof of the next lemma is the same as Lemma 2.2 of \cite{F}. 
\begin{lem}
\label{cex0}
The net $\big\{\sum_{z\in F}E_{z}AE_{z}\,\big|\,F\,\text{finite 
subset of}\,\s_{\mathop{\rm pp}}(U)\big\}$ 
converges in the strong operator 
topology.
\end{lem}
We symbolically write for such a limit
$$
\mathop{\rm s\!-\!lim}_{F\uparrow\s_{\mathop{\rm pp}}(U)}
\sum_{z\in F}E_{z}AE_{z}=:\sum_{z\in\s_{\mathop{\rm pp}}(U)}E_{z}AE_{z}\,.
$$
\begin{prop}
\label{procp}
Let $U$ be a unitary acting on the Hilbert space $\ch$, and 
$A\in\ck(\ch)$. Then
\begin{equation}
\label{wwlx}
\mathop{\rm w\!-\!lim}_{N\to+\infty}
\frac{1}{N}\sum_{n=0}^{N-1}
U^{n}AU^{-n}=\sum_{z\in\s_{\mathop{\rm pp}}(U)}E_{z}AE_{z}\,.
\end{equation}
\end{prop}
\begin{proof}
By approximating $A$ with a finite rank operator $A_{\eps}$, we have
\begin{align*}
&\bigg|\bigg\langle\frac{1}{N}\sum_{n=0}^{N-1}
U^{n}AU^{-n}x,y\bigg\rangle-\bigg\langle
\sum_{z\in\s_{\mathop{\rm pp}}(U)}E_{z}AE_{z}x,y\bigg\rangle\bigg|\\
\leq&\eps+\bigg|\bigg\langle\frac{1}{N}\sum_{n=0}^{N-1}
U^{n}A_{\eps}U^{-n}x,y\bigg\rangle-\bigg\langle
\sum_{z\in\s_{\mathop{\rm pp}}(U)}E_{z}A_{\eps}E_{z}x,y\bigg\rangle\bigg|\,.
\end{align*}

So, it is enough to check \eqref{wwlx} for rank one operators 
$A=\langle\,\cdot\,,\xi\rangle\eta$. In this situation, we have
\begin{align*}
\bigg\langle\frac{1}{N}&\sum_{n=0}^{N-1}U^{n}AU^{-n}x,y\bigg\rangle
=\frac{1}{N}\sum_{n=0}^{N-1}\langle U^{n}\eta,y\rangle
\langle U^{-n}x,\xi\rangle\\
&=\int_{\bt}\bigg(\frac{1}{N}\sum_{n=0}^{N-1}(z\bar w)^{n}\bigg)
\langle E(\di z)\eta,y\rangle\langle E(\di w)x,\xi\rangle\\
&\longrightarrow\bigg\langle
\sum_{z\in\s_{\mathop{\rm pp}}(U)}E_{z}AE_{z}x,y\bigg\rangle
\end{align*}
as 
$$\frac{1}{N}\sum_{n=0}^{N-1}(z\bar w)^{n}
\longrightarrow\chi_{1}(z\bar w) 
$$
pointwise. See Proposition 2.4 of \cite{F} for further details.
\end{proof}

\section{$C^{*}$--dynamical systems based on compact operators}

The present section is devoted to the study of some interesting 
ergodic properties of $C^{*}$--dynamical systems based on compact 
operators

Following the same lines of the previous results, we pass to the study 
of the convergence of Cesaro mean of automorphisms $\a$ of the 
$C^{*}$--algebra $\ck(\ch)$ consisting of all the compact operators 
acting on $\ch$. Consider the double transpose 
$\a^{**}\in\aut(\cb(\ch))$. As such an automorphism $\a^{**}$ is inner 
(cf. \cite{SZ}, Corollary 8.11), there exists a unitary $U$ acting 
on $\ch$ such that $\a=\ad{}\!_{U}$. Namely, each automorphism of 
$\ck(\ch)$ is implementable on $\ch$. 

Lemma \ref{cex0} allows us to define $\ce:\cb(\ch)\mapsto\cb(\ch)$ as 
\begin{equation}
\label{cex}
\ce(A):=\sum_{z\in\s_{\mathop{\rm pp}}(U)}E_{z}AE_{z}
\end{equation}
The properties of $\ce$ are collected in the following
\begin{prop}
\label{cexx}
The map $\ce$ is a conditional expectation projecting onto 
the $C^{*}$--subalgebra 
${\displaystyle\bigoplus_{z\in\s_{\mathop{\rm pp}}(U)}
E_{z}\cb(\ch)E_{z}}$.
\end{prop}
\begin{proof}
Following the same line of Lemma 2.1 of \cite{F}, we see that 
$\|\ce\|=1$. In addition
$$
\ce(\ce(A))=\sum_{z,w\in\s_{\mathop{\rm pp}}(U)}E_{w}E_{z}AE_{z}E_{w}
=\sum_{z\in\s_{\mathop{\rm pp}}(U)}E_{z}AE_{z}\equiv\ce(A)\,.
$$

Namely, $\ce$ is a norm one projection onto the 
the $C^{*}$--subalgebra 
${\displaystyle\bigoplus_{z\in\s_{\mathop{\rm pp}}(U)}
E_{z}\cb(\ch)E_{z}}$,
hence a conditional 
expectation, see \cite{St}, Theorem 9.1.
\end{proof}

Notice that the identity $\ce(I)$ of the range of $\ce$ is precisely
$$
E_{\mathop{\rm pp}}:=\sum_{z\in\s_{\mathop{\rm pp}}(U)}E_{z}\,,
$$
the selfadjoint projection onto the closed subspace of $\ch$ 
generated by the eigenvectors of $U$.

Now we specialize the matter to the case when $A$ is a compact 
operator.
\begin{lem}
\label{cecomp}
If $A\in\ck(\ch)$ then $\ce(A)\in\ck(\ch)$.
\end{lem}
\begin{proof}
We have by Schwarz, Holder and Bessel 
inequalities,
\begin{align*}
&|\langle\ce(A-B)x,y\rangle|\leq\sum_{z\in\s_{\mathop{\rm pp}}(U)}
|\langle(A-B)E_{z}x,E_{z}y\rangle|\\
\leq&\|A-B\|\sum_{z\in\s_{\mathop{\rm pp}}(U)}\|E_{z}x\|\|E_{z}y\|\\
\leq&\|A-B\|\big(\sum_{z\in\s_{\mathop{\rm pp}}(U)}
\|E_{z}x\|^{2}\big)^{1/2}
\big(\sum_{z\in\s_{\mathop{\rm pp}}(U)}
\|E_{z}y\|^{2}\big)^{1/2}\\
\leq&\|A-B\|\|x\|\|y\|\,.
\end{align*}

Thus, we can approximate $A$ by a finite rank operator. In addition, 
for a rank one operator $A=\langle\,\cdot\,,y\rangle x$, we have by 
polarization,
$$
A=\frac{1}{4}\sum_{\{z\in\bt\,:\,z^{4}=1\}}z\langle\,\cdot\,,x+zy\rangle(x+zy)\,.
$$

Namely, we can reduce the matter to the case when $A$ is the rank one 
positive operator $\langle\,\cdot\,,x\rangle x$. 
We now compute
$$
\langle\ce(A)x,x\rangle=\sum_{z\in\s_{\mathop{\rm pp}}(U)}
\|E_{z}x\|^{2}\,.
$$

As the last sum is convergent, there exists an at most countable set 
$z_{1},z_{2},\dots\subset\s_{\mathop{\rm pp}}(U)$ depending 
on $A$, such that $E_{z}x=0$ if $z\neq z_{j}$, $j=1,2,\dots\,\,$. In 
addition, ${\displaystyle\lim_{j}\|E_{z_{j}}x\|=0}$. Put 
$\l_{j}:=\|E_{z_{j}}x\|^{2}$ and 
$y_{j}:=\frac{E_{z_{j}}x}{\|E_{z_{j}}x\|}$, $j=1,2,\dots$. We have
$$
\ce(A)=\sum_{j}\langle\,\cdot\,,E_{z_{j}}x\rangle E_{z_{j}}x
=\sum_{j}\l_{j}\langle\,\cdot\,,y_{j}\rangle y_{j}\,.
$$

It readly seen that 
${\displaystyle\sum_{j=1}^{N}\l_{j}\langle\,\cdot\,,y_{j}\rangle 
y_{j}}$ converges in norm,
that is $\ce(A)$ is a compact operator.
\end{proof}
\begin{prop}
\label{cecomp1}
The restriction $E:=\ce\lceil_{\ck(\ch)}$ of the map in \eqref{cex} 
gives rise to a conditional expectation projecting onto the $C^{*}$--subalgebra 
$$
\bigg(\bigoplus_{z\in\s_{\mathop{\rm pp}}(U)}
E_{z}\cb(\ch)E_{z}\bigg)\bigcap\ck(\ch)\,.
$$
\end{prop}
\begin{proof}
Lemma \ref{cecomp} tells us that $\ce$ maps the compact operators into 
the compact ones. Moreover, $\big(\ce\lceil_{\ck(\ch)}\big)^{**}=\ce$ 
and the proof follows.
\end{proof}
\begin{thm}
\label{mmmmain}
Let $\a$ be an automorphism of $\ck(\ch)$, with $U$ the unitary acting 
on $\ch$ implementing $\a$. Then
$$
\lim_{N\to+\infty}\frac{1}{N}\sum_{n=0}^{N-1}\a^{n}=E\,,
$$
pointwise in the weak topology of $\ck(\ch)$, $E$ being the 
conditional expectation given in Proposition \ref{cecomp1}\,.
\end{thm}
\begin{proof}
By taking into account \eqref{wwlx},
\begin{equation}
\label{seq}
\mathop{\rm w\!-\!lim}_{N}\frac{1}{N}\sum_{n=0}^{N-1}\a^{n}(A)
=E(A)
\end{equation}
whenever $A$ is compact. Let now $T$ be a trace class operator and 
$T_{\eps}$ be a finite rank operator such that 
$\tr(|T-T_{\eps}|)\leq\eps$, $\tr$ being the unique normal faithtful 
semifinite trace on $\cb(\ch)$. Then
\begin{align*}
\tr\bigg(T\bigg(&\frac{1}{N}\sum_{n=0}^{N-1}\a^{n}(A)-E(A)\bigg)\bigg)
\leq\left|\tr\left((T-T_{\eps})\frac{1}{N}\sum_{n=0}^{N-1}\a^{n}(A)\right)\right|\\
+|\tr&((T-T_{\eps})E(A))|
+\tr\left(T_{\eps}\left(\frac{1}{N}\sum_{n=0}^{N-1}\a^{n}(A)-E(A)\right)\right)\\
&\leq2\eps\|A\|
+\tr\left(T_{\eps}\left(\frac{1}{N}\sum_{n=0}^{N-1}\a^{n}(A)-E(A)\right)\right)\,.
\end{align*}

Thus, we reduce the matter when $T$ is finite rank. The proof now follows 
by \eqref{seq}.
\end{proof}

Notice that if $\s_{\mathop{\rm pp}}(U)=\emptyset$, 
$$
\frac{1}{N}\sum_{n=0}^{N-1}\a^{n}\longrightarrow0\,,
$$
and if $\s_{\mathop{\rm pp}}(U)=\{1\}$ with $\Om$ the unique up to a 
phase invariant vector for $U$, that is in the case of weakly mixing $C^{*}$--dynamical 
systems based on compact operators, 
$$
\frac{1}{N}\sum_{n=0}^{N-1}\a^{n}\longrightarrow\om(\,\cdot\,)E_{1}\,,
$$
$\om$ being the vector state $\langle\cdot\,\Om,\Om\rangle$.

We end the present section with a recurrence result which is an immediate 
corollary of Proposition \ref{recc2}.

Let $(\ga,\a,\om)$ be a weakly mixing $C^{*}$--dynamical 
system, where $\ga=\ck(\ch)$, $\a=\ad{}\!_{U}$, and 
$\om=\langle\cdot\,\Om,\Om\rangle$, with $\Om$ invariant under the action of 
the unitary operator $U$. 
\begin{prop}
\label{rrreee}
Under the above conditions, if 
$\om(A)>0$, and $0<m_{1}<m_{2}<\cdots<m_{l}$ are 
natural numbers kept fixed, then there exists an $N_{0}$ such that 
$$
\frac{1}{N}\sum_{n=0}^{N-1}\om(A\a^{nm_{1}}(A)\a^{nm_{2}}(A)\cdots
\a^{nm_{l}}(A))>0
$$
for each $N>N_{0}$.
\end{prop}
\begin{proof}
By Proposition \ref{recc2}, we have
\begin{align*}
&\lim_{N\to+\infty}
\frac{1}{N}\sum_{n=0}^{N-1}
\om(A\a^{nm_{1}}(A)\a^{nm_{2}}(A)\cdots
\a^{nm_{l}}(A))\\
\equiv&\lim_{N\to+\infty}
\frac{1}{N}\sum_{n=0}^{N-1}\big\langle 
AU^{nm_{1}}AU^{n(m_{2}-m_{1})}\cdots A
U^{n(m_{l}-m_{l-1})}A\Om,\Om\big\rangle\\
=&\big\langle AE_{1}AE_{1}\cdots AE_{1}A\Om,\Om\big\rangle
\equiv\om(A)^{l+1}>0\,.
\end{align*}
\end{proof}

\end{document}